\documentstyle[11pt]{article}
\skip\footins 5mm plus 2mm minus 2mm
\begin{document}
\setlength{\parindent}{0 mm}
\newtheorem{theorem}{Theorem}
\newcommand {\athr} {\displaystyle \sum}
\vspace*{-0.5in}

\large

\centerline {{ \bf  MONIC POLYNOMIALS IN {\bf Z}[x] WITH ROOTS IN THE UNIT DISC   }} 

\bigskip
\centerline { {\bf Pantelis  A. Damianou}}

\vskip .5  cm

{\bf  A Theorem of Kronecker}

 This note
is motivated by an old result of Kronecker on  monic polynomials
with integer coefficients  having  all their roots in the unit disc.
We  call such polynomials {\it Kronecker polynomials} for short.
 Let $k(n)$ denote the number of Kronecker polynomials of degree $n$. We 
describe  a canonical form for such polynomials and use it to determine 
 the sequence $k(n)$,  for  small values of $n$.  The first step is to show that the number of Kronecker  polynomials 
of degree $n$ is finite. 
This fact is included in the following theorem  due to
Kronecker \cite{kronecker}. See also \cite{jones} for a more accessible
 proof.  The theorem actually gives more: the non-zero roots of such
polynomials are on the boundary of the  unit disc. We  use this fact
later on to show that these  polynomials  are essentially products of
cyclotomic polynomials.

\begin{theorem}
 Let $\lambda\not=0$ be a root of a monic polynomial $f(z)$ with integer coefficients. If all the roots 
of $f(z)$ are in the unit disc $\{ z \in {\bf C}\, | \, |z|\le 1 \}$,  then $|\lambda |=1$. 
\end{theorem}

\underline{Proof}

Let $n=$ deg$f$. The set of all monic polynomials of degree $n$   with integer coefficients having all their roots in the 
unit disc is finite.  To see this, we write
\begin{displaymath}
z^n+ a_{n-1}z^{n-1}+a_{n-2} z^{n-2}+ \cdots + a_0 = \prod_{j=1}^n \left( z-z_j \right) \ ,
\end{displaymath}
where $a_j\in {\bf Z}$ and $z_j$ are the roots of the polynomial.
 Using the fact that
$|z_j| \le 1 $ we have: 
\begin{displaymath}
\begin{array}{rcl}
|a_{n-1}|&=& |z_1+z_2+ \cdots +z_n| \le n = {n \choose 1}  \\
|a_{n-2}|&=& |\athr_{j,k} z_j z_k | \le    {n \choose 2} \\
 & \vdots &\\
|a_0| &=& |z_1 z_2 \dots z_n| \le 1= {n \choose n}    \ .
\end{array}
\end{displaymath}
Since the $a_j$ are integers, each $a_j$ is limited to at most $2   {n \choose j} +1$  values and therefore the 
total number of polynomials that satisfy the hypothesis of the theorem is finite.

We define 
\begin{displaymath}
f_1(z)=f(z)=\prod_{j=1}^n (z-z_j) \ ,
\end{displaymath}
and for $k\ge 2$ 
\begin{displaymath}
f_k(z)=\prod_{j=1}^n (z-z_j^k) \ .
\end{displaymath}
The fact that $f_k$ has integer coefficients follows from Newton's formula for symmetric polynomials:
If $s_k=\athr_{j=1}^n  z_j^k$, then
\begin{displaymath}
s_k-p_1 s_{k-1}+p_2 s_{k-2}- \cdots + (-1)^{k-1} p_{k-1} s_1 + (-1)^k k p_k =0   \ , 
\end{displaymath}
$k=1,2, \dots $. We denote by $p_k$  the elementary symmetric polynomial of degree $k$.
 It is also clear that all the roots of $f_k$ are in the unit disc. Therefore we  must have  $f_j(z)=f_k(z)$ for some $j <k$. Since the 
roots of $f_j$ are $z_1^j, z_2^j, \dots, z_n^j$ and the roots of $f_k$ are 
$z_1^k, z_2^k, \dots, z_n^k$ the two sets must be equal up to a permutation, i.e.,
\begin{displaymath}
z_r^j=z_{ \sigma (r)}^k  \ ,
\end{displaymath}
for $r=1,2,\dots, n $.  

Let $m$ be the order of $\sigma$ in $S_n$. Then we have for $z_1=\lambda$,
\begin{displaymath}
\lambda^{ j^m}=\lambda^{j j^{m-1}}=z_1^{j j^{m-1}}=z_{\sigma (1)}^{k  j^{m-1}}=
z_{\sigma (1)}^{j k j^{m-2}}=z_{ \sigma \left( \sigma(1) \right)}^{k^2 j^{m-2}}= \dots =\lambda^{k^m} \ .
\end{displaymath}
It follows that
\begin{displaymath}
\lambda^{j^m -k^m} =1 \ ,
\end{displaymath}
which implies that $\lambda $ is a root of unity.

{\bf Examples}

\begin{itemize}
\item
$ n=1$.  If $f(z)=z-a$, then $a$ takes the values $-1,0,1$. Therefore we have three linear 
Kronecker polynomials. They are $z$,
$z-1$, and  $z+1$. In this case $k(1)=3$.
\item
$ n=2$.  Let $f(z)=z^2-a z+b$.  We have $|a| \le 2$ and $|b| \le 1$. Therefore, $a=\pm 2, \pm 1, 0,$ and
$b=\pm 1, 0$. We have a total  of 15 polynomials, 9 of which are Kronecker polynomials. They are:
$z^2$, $z^2-z$, $z^2+z$, $z^2-1$, $z^2+z+1$, $z^2+1$, $z^2-z+1$, $z^2-2 z+1$, $z^2+2 z+1$. We obtain $k(2)=9$. 
\end{itemize}

This crude method is not very effective! For $n=3$, we have a total of $147$ candidates but only $19$ are
Kronecker polynomials. To make the counting more manageable we  need cyclotomic polynomials.

\bigskip

{\bf Cyclotomic polynomials}

It is easy to
see that the set of $n$th roots of unity  forms a cyclic group.  A complex number $\lambda$ is called a {\it primitive
$n$th root of unity} provided $\lambda$  is an $n$th root of unity and has order $n$ in the multiplicative group
of the $n$th roots of unity. In particular, a primitive root of unity generates the cyclic group of all $n$th
roots of unity. For example, if $n=3$ we have three cubic roots of unity: $1$, $\omega= { (-1 +\sqrt{3} i) \over 2} 
$ and  $\omega^{2}$; only  $\omega$ and $\omega^2$ are primitive cubic roots of unity. For $n=4$
we have two fourth primitive roots of unity, $\pm i$. 

We define the $n$th cyclotomic polynomial by
\begin{displaymath}
g_n(z)=(z-\lambda_1)(z-\lambda_2) \ldots (z-\lambda_s)  \ ,
\end{displaymath}
where $\lambda_1, \lambda_2, \cdots \lambda_s$ are all the distinct primitive $n$th roots of unity. For example, 
we easily calculate that $g_1(z)=z-1$, $g_2(z)=z+1$, $g_3(z)=z^2+z+1$, $g_4(z)=z^2+1$.  Since
 the $\lambda_j$ generate the cyclic group of roots of unity it follows that $s= {\rm deg}\,  g_n=\phi (n)$,  where $\phi$ is 
the Euler totient  function.

Let $s(j)$ denote the number of cyclotomic polynomials of degree $j$.  For each $j$, $s(j)$ is equal to
  the number of solutions
of $\phi (x)=j$.   In other words, it is the cardinality of the set $\phi^{-1} (j)$.
For example, $s(8)=5$ since the equation $\phi(x)=8$ has the following five solutions: 15, 16, 20, 24, and 30.  
We may have $s(j)=0$, for example, if $j$ is any odd integer bigger than $1$ but also for some even integers
 such as 14, 26, 34, 38, 50, \dots . Carmichael's conjecture is that $s(j) \not=1$ for all $j$; i.e.,  the equation
$\phi(x)=j$ has either no solution or at least two solutions.  Schlafly and
Wagon obtained  a lower bound  larger than $10^{10^7}$  for a counterexample to Carmichael's conjecture \cite{wagon}. 
An old conjecture of Sierpinski asserts that for each integer $k\ge 2$, there is an integer $j$ for which 
$s(j)=k$. This conjecture was proved recently by  Ford \cite{ford}.
  Erd\"os has shown   that
any value of the function $s$ appears infinitely often \cite{erdos}.

 Computing the value of $s(j)$ for a  particular $j$ (of 
moderate size) is not difficult; 
Suppose we want to calculate $s(2)$.  The first step is to find the solutions of the
equation $\phi(n)=2$. We write $n=\prod p_i^{\alpha_i}$ and, using the standard formula for $\phi(n)$,
 we obtain the equation

\begin{equation}
n \prod (p_i-1) = 2 \prod p_i \ .
\end{equation}
Let $d_i=p_i-1$. We search for integers $d_i$ that are divisors of $2$ and are such that   $d_i+1$ is
prime. The only solution is $d_1=1$, $d_2=2$. Therefore the possible  prime divisors of $n$ are $2$ and $3$.
We obtain  $n=2^{\alpha}3^{\beta}$, or $n=2^{\alpha}$, or $n=3^{\beta}$. 
In the first case we get  $2^{\alpha} 3^{\beta -1} =2$,  
which implies $\alpha=1$, $\beta=1$, i.e.,  $n=6$.  In the second case we have  $2^{\alpha -1} =2$, which
implies $\alpha=2$, i.e. $n=4$. In the last case we have $ 3^{\beta -1} 2 =2$, $\beta=1$ and $n=3$.
 Therefore $\phi^{-1} (2)=\{3, 4, 6 \}$, and   $s(2)=3$.  This is in principle the procedure one uses to calculate
$s(j)$ for small values of $j$. For more details see {\cite{burton}, p. 130}.
  With a symbolic manipulation package such as Maple one  may  easily compute   the values of 
$s$ up to billions (in fact, $ s(10^9)=152 $). 

 Returning to cyclotomic polynomials,  we record the important formula
\begin{equation}
z^n -1= \prod_{d|n} g_d(z) \ ,  
\end{equation}
where $d$ ranges over all positive divisors of $n$. 
This formula gives a recursive method for determining $g_n(z)$.  It is clear from (2)  that
whenever  $\lambda$ is a root of  unity, then $\lambda $ is also a root of some cyclotomic polynomial. The converse is
also true by definition.  We therefore obtain  the following:

\begin{theorem}
A  Kronecker polynomial is of the form
\begin{equation}
f(z)=z^k \prod_j g_j(z) \ ,  
\end{equation}
where $k \ge0$ is an integer and $g_j(z)$  are cyclotomic polynomials. 
\end{theorem}

In order to classify all Kronecker polynomials  of degree $n$, it  is enough to find all cyclotomic 
polynomials of degree $\le n$, and compute  all possible products of the form (3) that give a total degree
equal to $n$.  We illustrate with an example.

Let $n=3$. We need to  determine all  cyclotomic polynomials (together with $f(z)=z$) of degree less than or equal to $3$.
We have the  linear polynomials $f_1(z)=z$, $f_2(z)=g_1(z)=z-1$, and $f_3(z)=g_2(z)=z+1$. Since $s(2)=3$, we 
also have the following three quadratic  cyclotomic polynomials:   $h_1(z)=g_3(z)=z^2+z+1$, $h_2(z)=
g_4(z)=z^2+1$, and $h_3(z)=g_6(z)=z^2-z+1$.  There are no cubic cyclotomic polynomials since 
$s(3)=0$.  To obtain a Kronecker polynomial  of 
degree $3$, we have two possibilities:  Either it is a product of the form $f_i f_j f_k$ or a product of
the form $f_i h_k$. In the first case we have a total of ${ 5 \choose 3}=10$ choices. In the second case
we have a total of $3 \times 3=9$ choices. Therefore, the total number of Kronecker polynomials of degree
$3$  is $19$.  We obtain $k(3)=19$.

It seems  impossible to find a general, closed form formula for $k(n)$ since it requires an expression
for the multiplicity function $s$. Nevertheless, we  would like to outline  a procedure 
that can be  used  to calculate $k(n)$ for a specific  value of $n$.  The number of Kronecker polynomials of 
degree $n$ is the same as the number of solutions 
$\left( x_0,x_1, x_2, \dots  \right)$ to 
\begin{displaymath}
\sum_{r\ge 0} x_r \phi(r)=n  \ ,
\end{displaymath}
where for convenience, $\phi(0)=\phi(1)=1$. Since $\phi(2)=1$ and $\phi(r)$ is even for $r \ge 3$, we have 
the following expression for $k(n)$:
\begin{equation}
k(n)=\sum_{m=0}^n { n-m+2 \choose 2} b(m)   \ , 
\end{equation}
where $b(m)$ is the number of solutions to 
\begin{equation}
\sum_{r\ge 3} x_r \phi(r)=m  \ .  
\end{equation}
 The number of solutions of  (5)  for any particular $m$ (of moderate size) can be obtained   by  standard combinatorial formulas: One  
can write $b(m)$ in terms of certain  partitions of $m$ (consisting of even integers), and  the $s$--function as follows:
\begin{equation}
b(m)= \displaystyle{ \sum_{P} } \displaystyle{ \prod_j } { \mu(j,P)+s(j)-1 \choose \mu(j,P)} \ ,  
\end{equation}
where the sum is over all the partitions $P$ of $m$ with even parts, and $\mu(j,P)$ denotes the multiplicity
of $j$ as a part of $P$.  For example, $b(0)=1$, $b(2)=3$, $b(4)=10$, $b(6)=26$.

\smallskip
Using this  technique we  calculate that 
$k(1)=3$, $k(2)=9$, $k(3)=19$, $k(4)=43$, $k(5)=81$, $k(6)=159$, $k(7)=277$, $k(8)=501$, $k(9)=831$,
$k(10)=1415$, $k(11)=2253$, $k(12)=3673$, $k(13)=5675$, $k(14)=8933$, $k(15)=13447$, $k(16)=20581$,
$k(17)=30335$, $k(18)=45345$, $k(19)=65611$, $k(20)=96143$.

\bigskip

In closing we would like to mention another  approach  of calculating $k(n)$ based on the work of 
Boyd and Montgomery \cite{boyd}. Boyd and Montgomery determine the asymptotic
behaviour of the counting function for the number of Kronecker polynomials without zero as
a root. They also obtain a generating function for $k(n)$ of the form
\begin{equation}
\displaystyle{\prod_{n=0}^\infty} { 1 \over 1-x^{\phi(n)}}= \displaystyle{\sum_{n=0}^{\infty} k(n) x^n } \ .  
\end{equation}
One needs to know how fast $\phi(n)$ grows  in order to decide how many terms of the product to use. 
This estimate was  derived by Rosser and Schoenfeld in a well--known paper.  Using (7) and a computer program
it is easy to calculate $k(n)$ even for relatively large values of $n$. For example $k(100)=13445370780675$.

\smallskip
\noindent
{\bf Acknowledgments} I thank Stan Wagon and Herb Wilf for information on the
multiplicity function $s$; I thank Nickos  Papadatos for double--checking the values of
$k(n)$. Special thanks to the anonymous referee for  the beautiful formula (6).  I am grateful
to Cameron Stewart for pointing out reference \cite{boyd} and to  David  Boyd for explaining the contents
of \cite{boyd}.

\noindent
 {\it  University of Cyprus,  P. O. Box 20537,  1678 Nicosia, Cyprus} 

\noindent
{\it damianou@ucy.ac.cy}


\vskip 1cm
\begin{thebibliography}{99}

\bibitem{boyd}  D. W. Boyd, and H. L. Montgomery, {\it Number Theory, Banff, AB, 1988, Cyclotomic partitions},
de Gruyter, Berlin 1990, 7-25.

\bibitem{burton} D. M.  Burton, {\it Elementary Number Theory}, 3rd ed., McGraw Hill, New York, 1997.

\bibitem{erdos} P. Erd\"os, Some remarks on Euler's $\phi$-function, {\it Acta Arithm.} {\bf 4} (1958) 10-19.

\bibitem{ford} K. Ford, The number of solutions of $\phi(x)=m$, {\it Ann. of Math.} (2) {\bf 150} (1999)
283-311.


\bibitem{jones}  F. M. Goodman, P. de la Harpe, and  V. F. R. Jones, {\it Coxeter Graphs and Towers of Algebras},
Mathematical Sciences Research Institute Publications {\bf 14}, Springer-Verlag, New York, 1989.

\bibitem{kronecker}      L.  Kronecker,   Zwei S\"atze \"uber Gleichungen mit ganzzahligen Coefficienten, {\it Crelle},
Oeuvres I (1857) 105-108.

\bibitem{wagon} A. Schlafly, and S. Wagon, Carmichael's conjecture on the Euler function is valid below $10^{10000000}$,
{\it Math. Comp.} {\bf 63} (1994) 415-419.


\end{thebibliography}
\end{document}